\documentclass[reqno,11pt]{amsart}

\usepackage{amsmath,amssymb,mathtools,dsfont,graphicx}
\usepackage{graphicx}
\usepackage{mathptmx}

\usepackage{cite}

\title[Algorithm for numerical solutions to the kinetic equation of coalescing jumps]{
	Algorithm for numerical solutions to the kinetic equation of
	a spatial population dynamics model with coalescence and
	repulsive jumps
}

\author{
	Igor Omelyan
}
\address{Institute for Condensed Matter Physics, National Academy of Sciences of Ukraine, 
	1~Svientsitskii Street, UA-79011 Lviv, Ukraine}
\email{omelyan@icmp.lviv.ua}

\author{
	Yuri Kozitsky
}
\address{Instytut Matematyki, Uniwersytet Marii Curie-Skłodowskiej, 20-031 Lublin, Poland}
\email{jkozi@hektor.umcs.lublin.pl}

\author{
	Krzysztof Pilorz
}
\address{Instytut Matematyki, Uniwersytet Marii Curie-Skłodowskiej, 20-031 Lublin, Poland}
\email{krzysztof.pilorz@poczta.umcs.lublin.pl}

\begin{document}

\begin{abstract}
An algorithm is proposed for finding numerical solutions of a kinetic equation
that describes an infinite system of point articles placed in $\mathbb{R}^d (d \geq 1)$. 
The particles perform random jumps with pair wise repulsion, in the course of which they can 
also merge. The kinetic equation is an essentially nonlinear and nonlocal integro-differential 
equation, which can hardly be solved analytically. The derivation of the algorithm is based on 
the use of space-time discretization, boundary conditions, composite Simpson and trapezoidal 
rules, Runge-Kutta methods, adjustable system-size schemes, etc. The algorithm is then applied
to the one-dimensional version of the equation with various initial conditions. It is shown
that for special choices of the model parameters, the solutions may have unexpectable time
behaviour. A numerical error analysis of the obtained results is also carried out.
\end{abstract}
\maketitle

\section{Introduction}
Stochastic dynamics of infinite populations placed in a continuous habitat may include such
aspects as their motion accompanied by merging. A pioneering work in this direction was
published by Arratia in 1979, in which an infinite system of coalescing Brownian particles
in $\mathbb{R}$ was proposed and studied \cite{Arratia}. It was then extended to self-repelling 
motion of merging particles \cite{Toth}. Over the next two decades, different mathematical aspects 
of the Arratia model were intensively studied, see \cite{LeJan, Konarovskii, Berestycki, Renesse} and the 
references therein. In the Arratia flow, an infinite number of Brownian particles move 
independently up to their collision, then they coalesce and continue moving as single particles.

Recently, an alternative model of this kind has been proposed \cite{RCRJ1, RCRJ2}. Here, 
analogously to the Kawasaki approach \cite{Baranska, Berns}, particles make random jumps with 
repulsion acting on the target point. Additionally, two particles merge into one particle placed 
elsewhere with probability per time dependent on all the three locations and independent of the 
remaining particles. Thereafter, this new particle participates in the motion. The intensities and magnitude of jumps and coalescence are determined by the spatially nonlocal kernels. 

Similar models are used to describe predation in marine ecological systems, see, e.g.
Refs. \cite{Delius, Blanchard}. Identical processes are considered in freshwater ecosystems to 
study the phytoplankton dynamics \cite{Weisse, RudnickiWieczorek1, RudnickiWieczorek2}. Phytoplankton cells are 
dispersed in the water, leading to a patchy distribution of entities. The latter is the basis 
for the vast majority of oceanic and freshwater food chains. Another example concerns the modeling 
of cancer mechanisms according to one of which melanoma cells migrate and coalesce to form tumors 
\cite{Ambrose, Wessels}. Fresh tumor cells grow into clonal islands, or primary aggregates, that 
then coalesce to form heterogeneous formations. 

Quite recently \cite{RCRJ3}, first numerical results on the jump-coalescence kinetic equation 
have been reported. As a consequence, the role of repulsion interactions in the possible appearance 
of a spatial heterogeneity in the system has been elucidated. There it was shown also that the 
kinetic equation can rigorously be obtained from the corresponding microscopic evolution equations 
using previous theoretical approaches developed earlier 
\cite{MarkovEvolution, VlasovScaling, BanasiakLachowicz, SLM, PerturbativeApproaches, SPM} 
for other spatial population models. 
However, the numerical consideration was restricted to a few particular examples when choosing the 
parameters of the models and forms of the initial density. Moreover, very little was said about an 
algorithm used to obtain the numerical solutions.

In the present paper, we consistently derive the algorithm enabling to solve the 
repulsion-jump-coalescence kinetic equation. The derivation combines different techniques and 
develops methods allowing to reproduce properties of infinite systems on the basis of finite
samples. The algorithm is applied to population dynamics simulations of one-dimensional systems 
with various initial spatially inhomogeneous densities and forms of the jump, coalescence, and 
repulsion kernels. A comprehensive analysis of the obtained numerical results is also provided.

\section{Kinetic equation of a population model}
Consider an infinite population in the continuous space $\mathbb{R}^d$ of dimensionality 
$d \geq 1$. We will deal with two classes of stochastic events in the system: free coalescence 
and repulsive jumps. Time evolution of the population will be described in terms of the particle 
density $n(x,t)$. Performing the passage from the microscopic individual-based dynamics to the
mesoscopic description by means of a Vlasov-type scaling, one derives the following kinetic
equation for the density in the Poisson approximation with the jump-coalescence model 
\cite{RCRJ1, RCRJ2, RCRJ3}:
\begin{align}\label{eq:kineticGeneral}
	\frac{\partial n(x,t)}{\partial t} = &- \int a(x,y) \exp \left( - \int \varphi(y-u) n(u,t) du \right) n(x,t) dy \nonumber \\ 
	&+ \int a(x,y) \exp \left( - \int \varphi(x-u) n(u,t) du \right) n(y,t) dy \nonumber \\
	&- \int \int \Big[ b(x,y,z) + b(y,x,z) \Big] n(x,t) n(y,t) dy dz \nonumber \\
	&+ \int \int b(y,z,x) n(y,t) n(z,t) dy dz .
\end{align}

The first term in the rhs of Eq. \eqref{eq:kineticGeneral} represents the process of jumping of agents from
point $x$ to location $y$ with intensity $a(x, y)$ modulated by an exponential factor. This results
in a decrease of density $n(x,t)$ in point $x$ at time $t$. The exponential factor determines the
strength of repulsion between the particle jumping to $y$ and the rest of the system, where
$\varphi(y - u)$ denotes the corresponding interparticle potential. The intensity contribution and
repulsion strength depend on the configuration of all particles, so that the spatial integration
over $y$ and $u$ is carried out. The second term relates to the reverse process when particles
located anywhere in the coordinate space can appear at $x$ due to their random jumps to the
latter point. This increases the density $n(x,t)$. The third term defines the process of free 
coalescence when two particles, located at $x$ and $y$, merge into a single particle located at $z$
with intensities $b(x, y, z)$ and $b(y, x, z)$, resulting in a decrease of $n(x,t)$. Finally, merged particles can be randomly created in $x$ owing to the coalescence of arbitrary two other particles
located in $y$ and $z$ with intensity $b(y, z, x)$, leading to an increase of the local density.

The double integrations over $y$ and $z$ in the rhs of Eq. \eqref{eq:kineticGeneral} are performed in 
$\mathbb{R}^d$. They take into account the influence of all possible pair-wise coalescences on the 
change of $n(x,t)$ in $x$ at time $t$. These integrations are invariant with respect to the mutual 
replacement $y \leftrightarrow z$ in $b(x, y, z)$, $b(y, x, z)$, and $b(y, z, x)$. The jump and 
coalescence kernels are nonnegative functions symmetric w.r.t. first two arguments, 
$a(x, y) = a(y, x) \geq 0$ and $b(x, y, z) = b(y, x, z) \geq 0$. The kernels satisfy the 
integrability conditions $\int a(x, y)dx = \int a(x, y)dy = \mu_a$, $\int \varphi(x) = \mu_\varphi$, 
and $\int b(x_1, x_2, x_3 )dx_i dx_j = \mu_b$, where $i, j = 1, 2, 3$ with $i \neq j$. The quantities $\mu_a$, $\mu_b$, and $\mu_\varphi$ can be treated as parameters of the model. Moreover, according to
the translation invariance of the system, the following shifting identities should be satisfied,
$a(x, y) = a(x + u, y + u)$ and $b(x, y, z) = b(x + u, y + u, z + u)$. The obvious choice for the jump
kernel is $a(x, y) = a(x - y) = a(|x - y|)$.

The coalescence intensity will be modeled by the form
\begin{equation}\label{eq:coalescence}
	b(x,y,z) = b(x-y) \delta((x + y) / 2 - z),
\end{equation}
where $b(x - y) = b(|x - y|)$ and the Dirac $\delta$-function is applied. It implies that the 
resulting point of the coalescence of two particles in $x$ and $y$ is at the middle $z = (x + y)/2$ 
of their locations. This is quite natural for identical particles. Instead of $\delta$-function, 
we can, in principle, choose smoother dependencies, for example, Gaussian-like ones.

The $\delta$-function was used to simplify the calculations as then it allows us to lower 
dimensionality of the integration. Indeed, in view of Eq. \eqref{eq:coalescence} and integrating 
over $z$ to eliminate $\delta$-function, the kinetic equation \eqref{eq:kineticGeneral} transforms to
\begin{align}\label{eq:kinetic}
\frac{\partial n(x,t)}{\partial t} = &- \int a(x-y) \exp \left( - \int \varphi(y-u) n(u,t) du \right) n(x,t) dy \nonumber \\ 
&+ \int a(x-y) \exp \left( - \int \varphi(x-u) n(u,t) du \right) n(y,t) dy \nonumber \\
&- \int \left( 2 b(x-y) n(x,t) -2^d b(2(x-y)) n(2x - y, t) \right) n(y,t) dy,
\end{align}
where the symmetrical properties of the kernels have been taken into account. Imposing an
initial condition $n(x, 0)$, Eq. \eqref{eq:kinetic} leads to a complicated partial 
integro-differential equation with respect to the unknown function $n(x,t)$ at $t > 0$. 
Because of the presence of spatial integrals and nonlinearity, we doubt it can be solved 
analytically in general. Moreover, we cannot apply the convolution theorem to avoid the spatial 
integration in the last coalescence term of Eq. \eqref{eq:kinetic} due to the existence of three 
functions under the integrand which all depend on $y$. That is why we develop a numerical approach 
for this equation. The convolution method in the absence of coalescence is presented in Appendix. 
Analytical solutions in the spatially homogeneous case are also given there.

\section{Numerical algorithm}
\label{sec:Algorithm}
{\it Spatial discretization.} --- In order to solve numerically the kinetic equation it is necessary,
first of all, to perform its discretization in coordinate space. Let $x_i$ be the grid points uniformly distributed over $\mathbb{R}^d$ with mesh $h$ along all dimensions. For simplification of our
further presentation we will restrict ourselves in this study to a particular case of $d = 1$
(the generalization to any higher dimensionality $d > 1$ is straightforward and can be given
elsewhere). Then the discrete duplicate of Eq. \eqref{eq:kinetic} is
\begin{align}\label{eq:kineticDiscrete}
	\frac{dn_i}{dt} = &h \sum_j \Big( a_{i-j} \exp \Big[ -h \sum_k \varphi_{i-k}n_k \Big] n_j
	- a_{i-j} \exp \Big[ -h \sum_k \varphi_{j-k} n_k \Big] n_i \Big) \nonumber \\
	&- 2h \sum_j b_{i-j} n_i n_j + 2h \sum_j b_{2i - 2j} n_j n_{2i - j},
\end{align}
where $n_i \equiv n_i(t) = n(x_i,t)$ with $a_{i-j} = a(x_i - x_j) = a((i - j)h)$, 
$b_{i - j} = b(x_i - x_j) = b((i - j)h)$, $\varphi_{i-k} = \varphi(x_i - x_k) = \varphi((i - k)h)$, 
and the infinite sums over $j$ and $k$ represent the spatial integrals. It is obvious that in the 
limit $h \rightarrow 0$, the discretized kinetic equation \eqref{eq:kineticDiscrete} coincides
with its original continuous version \eqref{eq:kinetic}. Replacing $i - j$ by $j$, taking into account 
that the summation is carried out over the infinite number of terms, and introducing the auxiliary
quantities
\begin{equation}\label{eq:repulsionDiscrete}
	\lambda_i = \exp \Big[ -h \sum_k \varphi_{i-k} n_k \Big],
\end{equation}
one obtains that Eq. \eqref{eq:kineticDiscrete} can be cast more compactly in the following 
equivalent form
\begin{equation}\label{eq:kineticDiscreteCompact}
	\frac{dn_i}{dt} = h \sum_j \Big( a_j (\lambda_i n_{i-j} - \lambda_{i-j} n_i) - 2b_j n_i n_{i-j}
	+ 2b_{2j} n_{i-j} n_{i+j} \Big)
\end{equation}
where $a_j = a(jh)$, $b_j = b(jh)$, and $b_{2j} = b(2jh)$.

In computer simulations we cannot operate with infinite-size samples leading to the
infinite summation over j in Eqs. \eqref{eq:kineticDiscrete}, \eqref{eq:repulsionDiscrete}, 
and \eqref{eq:kineticDiscreteCompact}. Because of this we consider a finite
number $N$ of grid points $x_i$ uniformly distributed over the area $[-L/2, L/2]$ with spacing
$h = L/N$, where $i = 1, 2, \ldots, N$. This area will represent an interval, a square, or a cube in
cases $d = 1$, $2$, or $3$, respectively. The finite length $L$ should be sufficiently big with respect to
all characteristic coordinate scales of the system. The number $N$ of grid points must be large
enough to minimize the discretization noise. Then $h$ will be sufficiently small to provide a
high accuracy of the spatial integration. The finite-size effects can be reduced by applying
the corresponding boundary conditions (BC) when mapping infinite range $x \in ] - \infty, \infty [$ by
finite area $[-L/2, L/2]$. In view of the aforesaid, Eq. \eqref{eq:kineticDiscreteCompact} presents 
a coupled system of $N$ autonomous ordinary differential equations, where $i = 1, 2, \ldots, N$ and 
summation over $j$ is performed according to BC.

{\it Boundary conditions.} --- Three types of BC can be employed, namely, Dirichlet (DBC),
periodic (PBC), and asymptotic (ABC) boundary conditions. The choice depends on initial function 
$n(x, 0)$ and properties of solution $n(x,t)$. For example, if $n(x, 0)$ takes nonzero
values only within a narrow interval $[-l/2, l/2]$ with $l \ll L$, we can apply the DBC by letting 
$n_j = 0$ for all $|x_j| > L/2$. This means that during the finite simulation time $0 \leq t \leq T$,
the nonzero values of $n(x,t)$ do not approach the boundaries $x_B = \pm (L/2 - \max \sigma)$, where
$\max \sigma$ is the maximal radius of the kernels (see Sect. \ref{sec:Analysis}). In numerical 
calculations this can be expressed by the condition $n(x_B, t) < \varepsilon \ \max_x n(x,t)$, where 
$0 < \varepsilon \ll 1$ is the relative tolerable level (a negligibly small quantity slightly 
exceeding machine zero). When the propagation front becomes too close to the boundaries, i.e., 
$n(x_B, t) > \varepsilon \ \max_x n(x,t)$, we should enlarge $L$ (e.g. gradually doubling it) until 
to satisfy the required first condition, use DBC again, and continue the simulation for $t > T$. 
A special case, which is discussed in Sect. \ref{sec:Analysis}, where members of infinite configuration 
are initially absent in one half-space, requires a modified BC that combines DBC and ABC with 
an addition of adjustable system-size approach.

If $n(x, 0)$ and, thus, $n(x,t)$ are periodic functions, it is necessary to apply PBC with fixed
finite size $L$. According to PBC, the summation in Eq. \eqref{eq:kineticDiscreteCompact} for each 
current $i = 1, 2, \ldots, N$ is performed not only over all $j = 1, 2, \ldots, N$ but also over 
all infinite number of images $j'$ of $j$. The images are obtained by repeating the basic interval 
$[-L/2, L/2]$ by the infinite number of times to the left and to the right of it, so that 
$x_{j'} = x_j \pm KL$, where $K = 1, 2, \ldots, \infty$ and $n_{j'} = n_j$. This reproduces the 
periodicity $n(x \pm KL, 0) = n(x, 0)$, where $x \in [-L/2, L/2]$. The solution $n(x \pm KL,t) = n(x,t)$ 
will also be periodic for any time $t > 0$ with the same (finite) period $L$. In such a way the 
infinite system can be handled by a finite-size sample. Because the kernel values $a_j$ and $b_j$ 
decrease to zero with increasing the interparticle distance $|x_j|$, the summation over $j$ 
in Eq. \eqref{eq:kineticDiscreteCompact} can be actually restricted to a finite number of terms for
which $|x_j| \leq R_{a,b} < L/2$. The truncations radiuses $R_{a,b}$ are chosen to satisfy the conditions
$a(|x|) \approx 0$ and $b(|x|) \approx 0$ for $|x| > R_a$ and $|x| > R_b$, respectively.

In the spatially homogeneous case when $n(x,t) = n(t)$, we should apply ABC, i.e. $n_{j'} = n(t)$ 
for all $x_{j'} < -L/2$ and $n_{j'} = n(t)$ for all $x_{j'} > L/2$. For this case, PBC and ABC lead
to the same results. The ABC can also be used for spatially inhomogeneous solutions $n(x,t)$
which are flat for $x < - L/2$ and $x > L/2$ at a given $t$ where they take nonzero constant values.
Then $n_{j'} = n(-L/2,t)$ for all $x_{j'} < -L/2$ while $n_{j'} = n(L/2,t)$ for all $x_{j'} > L/2$. 
If in the course of time the flatness is violated at a current $L$, the basic length should be enlarged
using the automatically adjustable system-size approach mentioned above.

Taking into account the properties $a_{-j} = a_j$, $b_{-j} = b_j$ and the fact that the influence of
$a(x)$ and $b(x)$ can be neglected at large distances $x > R_{a,b}$, i.e. assuming that $a(x) = 0$ for
$|x| > R_a$ and $b(x) = 0$ for $|x| > R_b$, Eq. \eqref{eq:kineticDiscreteCompact} transforms to
\begin{align}\label{eq:discretized}
	\frac{dn_i}{dt} = &h \sum_{j=1}^{j_a} \xi_j^{(a)} a_j \Big( \lambda_i (n_{i-j} + n_{i+j}) - 
	n_i (\lambda_{i-j} + \lambda_{i+j}) \Big) \nonumber \\
	&- 2h \xi_0^{(b)} b_0 n_i^2 -2h n_i \sum_{j=1}^{j_b} \xi_j^{(b)} b_j (n_{i-j} + n_{i+j}) \nonumber \\
	&+ 2h \xi_0^{(2b)} b_0 n_i^2 +4h \sum_{j=1}^{j_b/2} \xi_j^{(2b)} b_{2j} n_{i-j} n_{i+j}
\end{align}
where $i = 1, 2, \ldots, N$ and the summations over $j$ are performed already with finite positive
nonzero integers up to $j_a = R_a / h$ or $j_b = R_b / h$. Quite similarly, using the properties 
$\varphi_{-j} = \varphi_j$ and $\varphi(x) = 0$ for $|x| > R_\varphi$, one finds from Eq. \eqref{eq:repulsionDiscrete} that
\begin{equation}\label{eq:discretizedRepulsion}
	\lambda_i = \exp \Big[ -h \sum_j^{BC} \xi_j^{(\varphi)} \varphi_j n_{i-j} \Big] = \exp 
	\Big[ -h \xi_0^{(\varphi)} \varphi_0 n_i - h \sum_{j=1}^{j_\varphi} \xi_j^{(\varphi)} \varphi_j
	(n_{i-j} + n_{i+j}) \Big],
\end{equation}
where $i = 1, 2, \ldots, N$ and $j_\varphi = R_\varphi /h$ with truncation radius $R\varphi < L/2$. 
It is evident that in the limits $L, N, R \rightarrow \infty$ provided $h \rightarrow 0$, the 
discretized equation \eqref{eq:discretized} with Eq. \eqref{eq:discretizedRepulsion} coincide with its 
original, continuous counterpart \eqref{eq:kinetic}.

Weights $\xi_j^{(a,b,2b,\varphi)}$ at kernel values $a_j, b_j, b_{2j},$ and $\varphi_j$ were introduced 
in Eqs. \eqref{eq:discretized} and \eqref{eq:discretizedRepulsion} to improve precision of the numerical integration over coordinate space. They are determined according to the chosen method \cite{NumericalMethods}. 
These weights satisfy the normalization condition $\xi_0 + 2 \sum_{j=1}^{j^*} \xi_j = 2 j^*$, where
$j^* = j_{a,b,\varphi}$ or $j^* = j_b /2$. In particular, we can use the composite Simpson or 
trapezoidal rules. For example, in the composite trapezoidal scheme
of the second order we have that $\xi_0 = \xi_1 = \xi_2 = \ldots = \xi_{j^*-1} = 1$, while 
$\xi_{j^*} = 1/2$. The composite Simpson method of the fourth order yields 
$\xi_{j^*} = 1/3, \xi_{j^* - 1} = 4/3, \xi_{j^* - 2} = 2/3, \xi_{j^* - 3} = 4/3, \xi_{j^* - 4} = 2/3, 
\ldots, \xi_0 = 4/3$ if $j^*$ is odd and $\xi_0 = 2/3$ if $j^*$ is even. The latter is
more accurate than the former and involve numerical uncertainties of order of $\mathcal{O}(h^4)$ 
versus $\mathcal{O}(h^2)$.

We mention that $n_i$ are explicitly defined at $i = 1, 2, \ldots, N$, so that the corresponding
BC should be applied in Eqs. \eqref{eq:discretized} and \eqref{eq:discretizedRepulsion} to $n_{i-j}$ 
and/or $n_{i+j}$ whenever $i - j < 1$ and/or $i + j > N$, because $j = 1, 2, \ldots, j^* < N/2$. 
The uniform knot distribution over $[-L/2, L/2]$ can be chosen in the form $x_i = -L/2 + (i - 1/2)h$, 
where $i = 1, 2, \ldots, N$ with even $N$. This provides the symmetricity of knot positions with 
respect to $x = 0$. Then according to PBC the calculations of $n_{i \pm j}$ should be performed as
\begin{equation}\label{eq:PBC}
	n_{i-j} = \left\{
	\begin{array}{lr}
		n_{i-j+N}, & i - j < 1\\
		n_{i-j}, & 1 \leq i - j \leq N
	\end{array}
	\right. ,\quad n_{i+j} = \left\{
	\begin{array}{lr}
		n_{i+j-N} & i + j > N\\
		n_{i+j} & 1 \leq i - j \leq N
	\end{array}
	\right. .
\end{equation}
Note that $i = 1, 2, \ldots, N$ and $j = 1, 2, \ldots, j^* < N/2$. The applications of DBC and ABC result
in
\begin{equation}\label{eq:DBC}
n_{i-j} = \left\{
\begin{array}{lr}
0, & i - j < 1\\
n_{i-j}, & 1 \leq i - j \leq N
\end{array}
\right. ,\quad n_{i+j} = \left\{
\begin{array}{lr}
0 & i + j > N\\
n_{i+j} & 1 \leq i - j \leq N
\end{array}
\right. .
\end{equation}
and
\begin{equation}\label{eq:ABC}
n_{i-j} = \left\{
\begin{array}{lr}
n_1, & i - j < 1\\
n_{i-j}, & 1 \leq i - j \leq N
\end{array}
\right. ,\quad n_{i+j} = \left\{
\begin{array}{lr}
n_N & i + j > N\\
n_{i+j} & 1 \leq i - j \leq N
\end{array}
\right. .
\end{equation}
respectively. In the cases $d = 2$ and $d = 3$, the boundary transformations can be implemented
similarly to those given by Eqs. \eqref{eq:PBC} -- \eqref{eq:ABC}.

{\it Time integration.} --- In the most general form, our coupled system of $N$ autonomous
ordinary differential equations can be given as
\begin{equation}\label{eq:discretizedGeneral}
	\frac{d n_i}{d t} = \dot{n}_i = s_i(n_1, n_2, \ldots n_N),
\end{equation}
where $i = 1, 2, \ldots, N$ and $s_i (n_1, n_2, \ldots, n_N )$ represents the rhs of Eq. 
\eqref{eq:discretized} with taking into account Eq. \eqref{eq:discretizedRepulsion}. 
Introducing the set $\Gamma(t) = \{n_i (t)\}$ of $N$ dynamical variables and their time
derivatives $\Phi(t) = \{s_i (t)\}$, Eq. \eqref{eq:discretizedGeneral} can be compactly cast as 
$\dot{\Gamma} = \Phi(\Gamma)$. A common practice to solve differential equations of such a type 
is to use classical Runge-Kutta schemes \cite{NumericalMethods}. The scheme of the fourth 
order (RK4) reads
\begin{equation}\label{eq:RK4}
	\Gamma(t + \Delta t) = \Gamma(t) + \frac{\Delta t}{6} \Big( \Phi_1 + 2 \Phi_2 + 2 \Phi_3 + \Phi_4)
	+ \mathcal{O}(\Delta t^5),
\end{equation}
where $\Delta t$ is the time step and 
$\Phi_1 = \Phi(\Gamma (t))$, $\Phi_2 = \Phi(\Gamma (t) + \Phi_1 \Delta t/2)$, 
$\Phi_3 = \Phi(\Gamma(t) + \Phi_2 \Delta t/2)$, $\Phi_4 = \Phi(\Gamma(t) + \Phi_3 \Delta t)$ are the 
derivatives in intermediate stages. The Runge-Kutta scheme of the second order (RK2) can also be applied. 
There are two versions of RK2. The first one (will be referred to as simply RK2) is known as the Heun 
method and can be related to the trapezoidal or Verlet-like integration. It has the form 
$\Gamma (t + \Delta t) = \Gamma (t) + [\Phi(\Gamma (t)) + \Phi(\Gamma (t) + \Phi (\Gamma(t)) \Delta t)] 
\Delta t/2 + \mathcal{O}(\Delta t^3)$. The second version (RK2') relates to a
middle point scheme, $\Gamma (t + \Delta t) = \Gamma (t) + \Phi(\Gamma (t) + \Phi(\Gamma (t)) \Delta t/2)
\Delta t + \mathcal{O}(\Delta t^3)$. The RK2 and RK20 approaches are two-stage schemes. 
The RK4 integration consists of four stages, so that at the same $\Delta t$ it will require twice larger 
number of operations to cover the same time interval $T$ of the observation with respect to RK2. 
However, RK4 is more accurate than RK2 with $\mathcal{O}(\Delta t^5)$ versus $\mathcal{O}(\Delta t^3)$ 
local truncation errors.

In such a way, the numerical solution $n_i (t)$ can be found for any $t = p \Delta t$, where 
$p = 1, 2, \ldots, P$ over the interval $t \in [0, T]$ with $T = P \Delta t$ by sequentially 
repeating $P$ times the transformations given by Eq. \eqref{eq:RK4}. The 
$\mathcal{O}(\Delta t^5 )$-uncertainties can be reduced to an arbitrary
small level by decreasing the length $\Delta t$ of the time step. We mention that the finite-size
effects are minimized by choosing a sufficiently large size $L$ of the system and applying
the corresponding boundary conditions. The uncertainties caused by the discretization are
reduced by decreasing mesh $h = L/N$. The latter can be achieved at sufficiently large values
of $N \gg 1$. In general, the total number of operations per one time step, which are necessary
to obtain solutions $n_i (t)$ for the set of $N$ differential equations \eqref{eq:discretizedGeneral}, 
is proportional to $N^2$. For $d > 1$, the computational efforts will increase drastically with 
increasing $d$, namely as $N^{2d}$. For convolution solutions this number is lowered to order of 
$N^d \ln N^d$ (see Appendix). In the case of simple rectangle kernels the total number of operations 
can be decreased to be proportional to $N^d$.

\section{Applications and analysis of the results}
\label{sec:Analysis}
{\it Numerical details.} --- The numerical simulations were carried out for one-dimensional systems ($d = 1$) of free or repulsive jumping particles which can coalesce. The kinetic equation
\eqref{eq:kinetic} was solved by using the algorithm described in the preceding section. The discretization
in $\mathbb{R}$ was done with a mesh of $h = 0.0125$ -- $0.2$ in dependence on the choice of initial conditions and kernel parameters. The initial basic length was $L = 20$ within either the fixed (for PBC) or adjustable-size (for DBC and ABC) regime. Spatial integration was performed
by employing the composite Simpson rule. Time integration was done with the help of the
RK4 scheme at a step of $\Delta t = 0.1$. The relative tolerable level was chosen to be equal to
$\varepsilon \sim 10^{-12}$. Further increasing space and time resolutions did not noticeably affect the solutions (see the last subsection).

The jump $a(x)$, coalescence $b(x)$, and repulsion $\varphi(x)$ kernels were modelled by Gaussian
\begin{equation} \label{eq:Gaussian}
	G_{\mu, \sigma} (x) = \frac{\mu}{(2 \pi \sigma^2)^\frac{1}{2}} \exp \left( - \frac{x^2}{2 \sigma^2} \right)
\end{equation}
or rectangle
\begin{equation} \label{eq:Rectangle}
	C_{\mu, \sigma} (x) = \left\{
	\begin{array}{rcl}
		\frac{\mu}{2\sigma}, & & |x| \leq \sigma\\
		0, & & |x| > \sigma
	\end{array}
	\right.
\end{equation}
functions, where $\mu = \mu_a$, $\mu_b$, or $\mu_\varphi$ and $\sigma = \sigma_a$, $\sigma_b$, or 
$\sigma_\varphi$ are the intensities and ranges of the corresponding interactions, respectively. 
The kernels are normalized, $\int \{a, b, \varphi\}(x)dx = \mu_a, \mu_b, \mu_\varphi$. A symmetrical pair of shifted Gaussian or rectangle functions was involved as well,
\begin{equation} \label{eq:Shift}
	F_{\mu, \sigma, s}(x) = \frac{1}{2} \Big( F_{\mu, \sigma}(x-s) + F_{\mu, \sigma}(x+s) \Big),
\end{equation}
where $F \equiv G$ or $C$ and $s$ is the shifting interval. The truncation radius for Gaussian kernels
was $R = Q \sigma$ with $Q = 6$ for which $G_{\mu,\sigma} (Q \sigma)/G_{\mu,\sigma} (0) \sim 10^{-8}$, 
so that their contribution at $|x| > R$ can be neglected. In the case of discontinuous rectangle kernels 
we have $Q = 1$ by definition. For the sums of two shifted kernels \eqref{eq:Shift} the truncation 
distance increases to $R = Q \sigma + s$.

In view of the above, we have six kernel parameters, $\mu_a$, $\mu_b$, and $\mu_\varphi$, as well as 
$\sigma_a$, $\sigma_b$, and $\sigma_\varphi$ for the repulsion-jump-coalescence model. 
This leads to a large number of all possible combinations. 
Because of this we will deal with most characteristic examples when
choosing values for these parameters for single Gaussian and rectangle kernels or their $\pm s$-shifted 
double counterparts. In order to consistently analyze their influence on the dynamics
of populations the following four situations will be considered: (i) pure free jumps; (ii)
repulsive jumps; (iii) pure coalescence; and (iv) repulsive jumps with coalescence. It is worth
emphasizing also that solution $n(x,t)$ depends cardinally on the choice of initial condition
$n_0 (x) \equiv n(x, 0)$. Like kernels, the initial condition can be selected in the form of Gaussians
\eqref{eq:Gaussian} or rectangles \eqref{eq:Rectangle}, trigonometric or step functions, etc.

\begin{figure}
	\centering
	\includegraphics{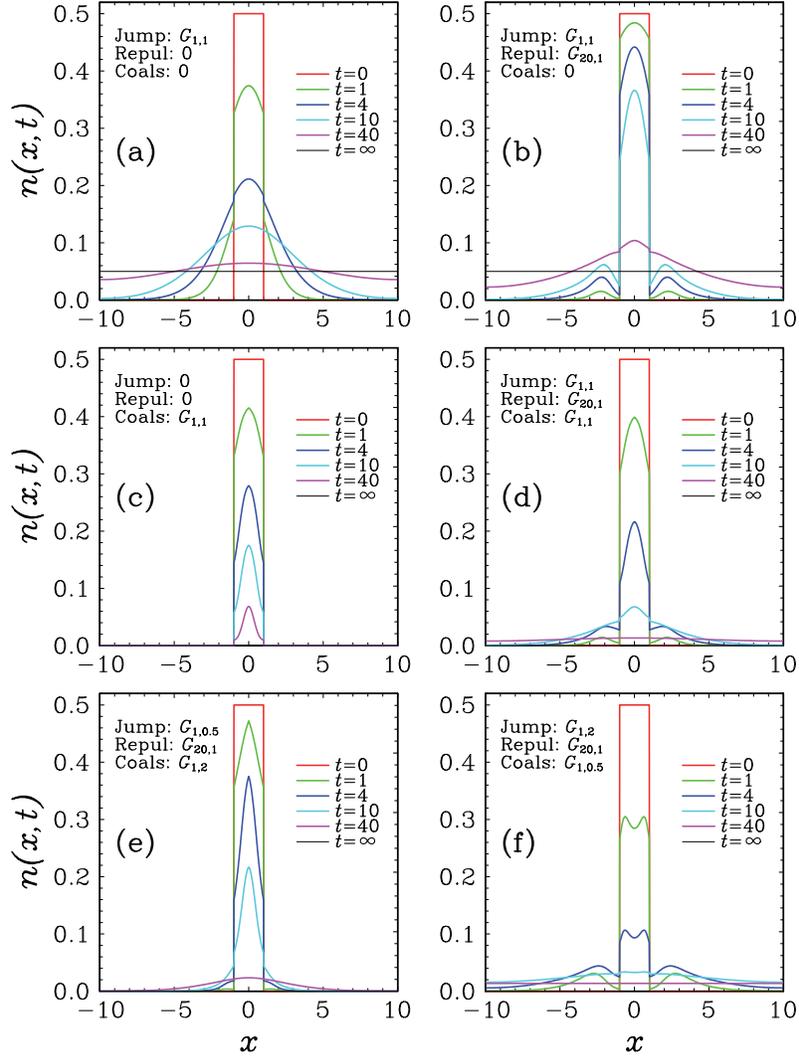}
	\caption{
		Time evolution of density $n(x,t)$ as dependent on spatial coordinate $x$ at several 
		moments of time $t$ for	rectangle initial condition $n(x, 0) = C_{1,1,L}(x)$ in the cases: 
		(a) pure free jumps; (b) repulsive jumps; and (c) pure coalescence; as well as repulsive 
		jumps and coalescence for (d) equal interaction ranges; (e) short-ranged jumps; 
		and (f) short-ranged coalescence. The infinite system is periodic ($L = 20$) on the interval 
		$x \in [-10, 10]$ with jump, repulsion, and coalescence Gaussian kernels of different intensities 
		and ranges (see the legends	inside).
	}
	\label{fig:1}
\end{figure}

{\it Rectangle initial density profiles.} --- The first example relates to the initial condition in
the form of periodic rectangle function $n(x, 0) = C_{1,1,L} (x)$. The infinite system is reproduced
by repeating the single rectangle segment $C_{1,1}(x)$ at $x \in \Omega = [-L/2, L/2]$ on the interval
$]-\infty, \infty[$ with a period of $L$ and applying PBC. The jump $a(x) = G_{\mu_a, \sigma_a}(x)$, 
repulsion $\varphi(x) = G_{\mu_\varphi,\sigma_\varphi} (x)$, and coalescence $b(x) = G_{\mu_b, \sigma_b} (x)$ 
kernels are modelled by the Gaussians with intensities $\mu_a = 1$, $\mu_\varphi = 20$, and $\mu_b = 1$. 
Three sets of kernel ranges were considered, namely,
$\sigma_a = \sigma_\varphi = \sigma_b = 1$, $\sigma_a < \sigma_\varphi < \sigma_b$ with $\sigma_a = 0.5$, 
$\sigma_\varphi = 1$, $\sigma_b = 2$ and vice versa, $\sigma_a > \sigma_\varphi > \sigma_b$ with 
$\sigma_a = 2$, $\sigma_\varphi = 1$, $\sigma_b = 0.5$. The corresponding time evolution of spatial 
structure $n(x,t)$ is presented in Fig. \ref{fig:1} for the cases of pure free jumps, repulsive jumps, 
pure coalescence, and repulsive jumps with coalescence with $\sigma_a = \sigma_\varphi = \sigma_b = 1$, 
parts (a), (b), (c), and (d), respectively, as well, with $\sigma_a < \sigma_\varphi < \sigma_b$ and 
$\sigma_a > \sigma_\varphi > \sigma_b$ for coalescing repulsive jumps, parts (e) and (f).

From Fig. \ref{fig:1}(a) we see that for pure free jumps, the discontinuity of $C_{1,1}(x)$ at 
$x = \pm 1$ soon disappears with increasing time, transforming the initial density into a 
Gaussian-like shape at $t \gtrsim 10$. For longer times, $t \gtrsim 40$, the solution $n(x,t)$ 
tends to more and more homogeneous densities. In the limit $t \rightarrow \infty$ we expect 
absolutely flat function $\lim_{t \rightarrow \infty} n(x,t) = n$ independent on $x$, where
$n = \frac{1}{L} \int_{\Omega} n(x, 0) dx = 1/L = 0.05$. Moreover, for any time $t$ the total 
number $\int_{\Omega} n(x,t)dt$ of particles on the interval $\Omega = [-L/2, L/2]$ is 
constant because of the absence of coalescence. For strong repulsive jumps [part (b)] 
the density function $n(x,t)$ remains discontinuous at $x = \pm 1$ up to $t \sim 40$ and 
its shape at shorter times is more complicated. In particular, apart from the existence 
of the main maximum at $x = 0$ which decreases in amplitude with increasing $t$, 
two symmetric secondary maxima appear additionally at $0 < t \lesssim 40$ in the ranges 
$x \approx \pm 2$ due to the repulsion between particles. We believe that all the maxima 
disappear at $t \rightarrow \infty$ with the same asymptotic behaviour 
$\lim_{t \rightarrow \infty} n(x,t) = n = 1/L = 0.05$ as for free jumps. In contrast, 
for free coalescence [part (c)] we expect decay of $n(x,t)$ to zero at $t \rightarrow \infty$. 
Moreover, here the particles remain to be located exclusively within the initial interval $[-1, 1]$ 
and they are absent outside of it at any $t$. In other words, no density propagation is indicated 
because the particles do not move apart from coalescence.

When the repulsive jumps are carried out in the presence of coalescence at equal interaction ranges 
$\sigma_a = \sigma_\varphi = \sigma_b = 1$, see part (d) of Fig. \ref{fig:1}, the pattern is somewhat 
similar to that of part (b). However, the three-maximum structure dissipates now much faster, leading to 
the zeroth asymptotics already at $t \gtrsim 40$. For short-ranged jumps, where $\sigma_a = 0.5$,
$\sigma_\varphi = 1, \sigma_b = 2$, the central peaks at $x = 0$ become sharper, while the secondary 
side maxima at $x \approx \pm 2$ do not appear, see part (e) and compare it with subsect (d). 
In the case of short-ranged coalescence when $\sigma_a = 2, \sigma_\varphi = 1, \sigma_b = 0.5$, 
the central peaks transform into a more complicated structure with one central minimum at $x = 0$ and 
two side maxima at $x \approx \pm 0.5$, see part (f). The secondary maxima at $x \approx \pm 2$ becomes 
more visible with respect to those for equal-range interactions [cf. part (f)]. Thus, the influence 
of jumps on the dynamics increases not only with increasing their intensity but range as well. The same
concerns the coalescence. Note also that the density profiles in Figs. \ref{fig:1} (a)--(f) are symmetric,
i.e., $n(-x,t) = n(x,t)$, like the initial condition, $n(-x, 0) = n(x, 0)$. This follows from the
symmetry of the kinetic equation.

\begin{figure}
	\centering
	\includegraphics{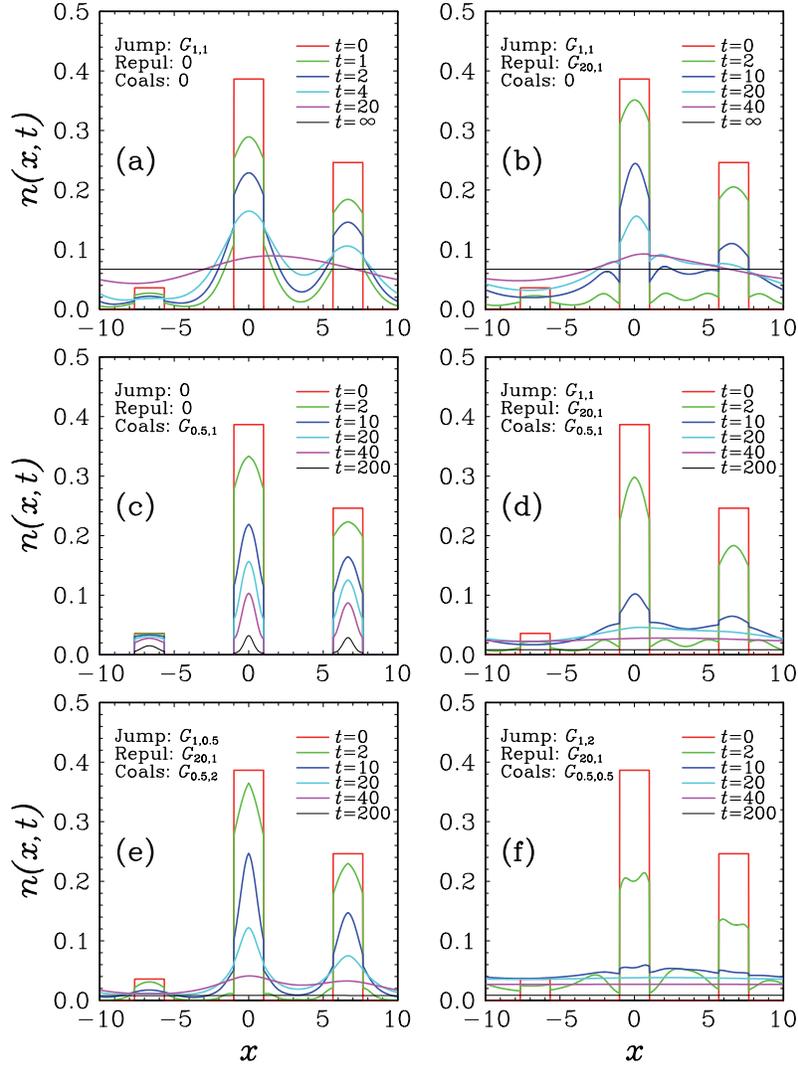}
	\caption{
		Time evolution of solution $n(x,t)$ for asymmetric initial density in the form of three 
		rectangle functions with different amplitudes. Other notations are the same as in 
		Fig. \ref{fig:1}.
	}
	\label{fig:2}
\end{figure}

The second choice deals with asymmetric initial condition $n(-x, 0) \neq n(x, 0)$ in the form
of $\mathcal{N}_0$ shifted single rectangle functions $C_{v_k, \sigma_k} (x + s_k)$ with intensities 
$v_k$ and ranges $\sigma_k$, namely,
\begin{equation}\label{eq:MultiShift}
	n(x,0) = C_{v_k,\sigma_k, s, \mathcal{N}_0} (x) = \sum_{k=1}^{\mathcal{N}_0} C_{v_k, \sigma_k} (x + s_k)
\end{equation}
where $s_k = -L/2 + (k - 1/2)L/\mathcal{N}_0$ are shifting parameters. Repeating \eqref{eq:MultiShift} 
with period $L$, we should apply PBC to deal with the infinite system. We used a particular case of 
\eqref{eq:MultiShift} with $\mathcal{N}_0 = 3$ and $L = 20$ as well as three different amplitudes 
$v_{1,2,3}$ generated at random in the interval $]0, 1[$. The corresponding result on this is shown in 
Fig. \ref{fig:2}.

Looking at Fig. \ref{fig:2} and comparing it with Fig. \ref{fig:1} we see that behaviour of $n(x,t)$ 
at short times can be approximately presented as a sum of independent separate solutions obtained
for single rectangle initial densities $C_{v_k, \sigma_k} (x + s_k)$. With increasing time, 
a coherence between the separate solutions appears. Again, in the absence of coalescence [$\mu_b = 0$, 
parts (a) and (b)] the density $n(x,t)$ at $t \rightarrow \infty$ flattens in $x$ and seems tend to 
the nonzero constant $n = \frac{1}{L} \int_{\Omega} n(x, 0) dx = (v_1 + v_2 + v_3)/L$ for any 
initial distributions. We mention that at $\mu_b = 0$, the total amount $\int_{\Omega} n(x,t)dx$ 
of particles in $\Omega$ is unchanged and equal to its initial value $\int_{\Omega} n(x, 0)dx$. 
For $\mu_b > 0$, the density function $n(x,t)$ seems to take its zeroth asymptotics 
$\lim_{t \rightarrow \infty} n(x,t) = 0$ at $t \rightarrow \infty$ [parts (c)--(f)] with the flow 
of time (except special choices, see below). For populations with an infinite number of particles 
and periodic initial conditions $n(x \pm KL, 0) = n(x, 0)$ with $x \in [-L/2, L/2]$, where 
$K = 1, 2, \ldots, \infty$, the solution $n(x \pm KL,t) = n(x,t)$ will also be periodic for any 
time $t > 0$ with the same (finite) period $L$. Then, in particular, $n(-L/2,t) = n(L/2,t)$ as 
is confirmed in Fig. \ref{fig:2}. Investigations show that the increase of the strength $\mu$ 
and range $\sigma$ of the jump and coalescence kernels accelerates this process, see for instance, 
parts (e) and (f) of Figs. \ref{fig:1} and \ref{fig:2}. Using the Gaussian (instead
of rectangle) initial conditions and rectangle (instead of Gaussian) kernels leads to results
(not shown) which are similar to those of Figs. \ref{fig:1} and \ref{fig:2}.

\begin{figure}
	\centering
	\includegraphics{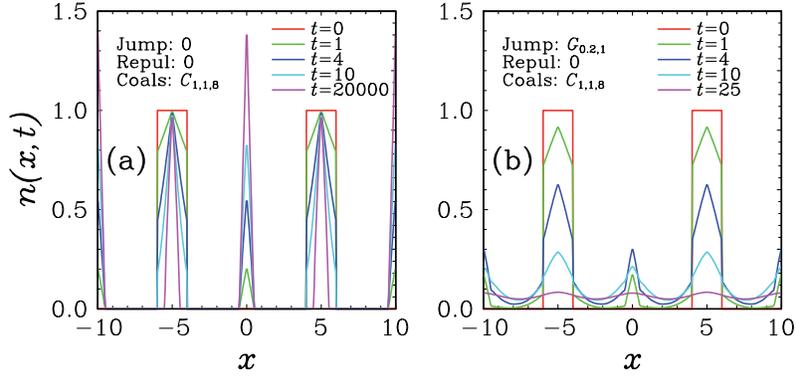}
	\caption{
		Dynamics of spatial structure $n(x,t)$ for initial density in the form of two rectangle functions,
		$n(x, 0) = C_{1,1,5,2}(x)$ in the cases: (a) pure coalescence and (b) coalescence with free jumps. 
		The interactions are modelled by shifted pair rectangle coalescence $C_{1,1,8}(x)$ and 
		Gaussian jump $G_{0.2,1}(x)$ kernels.
	}
	\label{fig:3}
\end{figure}

As visualized in Figs. \ref{fig:1} and \ref{fig:2}, the presence of coalescence ($b(x) \neq 0$) 
seems to lead to the zeroth asymptotic $\int_{\Omega} n(x,t)dx \rightarrow 0$ at long times 
provided the kernels are single rectangle functions with positive values around zero 
(the same concerns simple Gaussians). However, the coalescence kernel can be chosen in the form 
$b(x) = C_{\mu,\sigma,s}(x)$ of a pair of two shifted rectangle functions [Eq. \eqref{eq:Shift}] 
with appropriate shifting parameter $s$ to avoid the zeroth density limit. The initial inhomogeneous 
density $n(x, 0)$ over the basic interval $[-L/2, L/2]$ with PBC at period $L$ should also be chosen 
correspondingly. For instance, we can consider $n(x, 0)$ in the form $C_{v_k,\sigma_k,s',\mathcal{N}_0} (x)$ 
of two ($\mathcal{N}_0 = 2$) single rectangle functions [Eq. \eqref{eq:MultiShift}] with the same amplitude 
$v_k = \mu_b = 1$ and ranges $\sigma_k = \sigma_b \equiv \sigma = 1$ as those of the coalescence kernel 
$b(x)$ but with different shifting parameter $s' \neq s$ satisfying the constraint $s' = 2s - 2 \sigma$. 
Choosing $s = 5$ one finds at $\sigma = 1$ that $s' = 8$. The corresponding result
for the case $n(x, 0) = C_{1,1,5,2} (x)$ with $L = 20$ and $b(x) = C_{1,1,8} (x)$ in the absence 
($a = 0$) of jumps is depicted in Fig. \ref{fig:3}(a). We see that at the beginning, the rectangles soon transform into triangle-shaped peaks centered at $x = \pm (5 + KL)$, while the additional triangles appear
exactly in the middle of them at $x = 0 \pm KL$ (below we will omit the terms $\pm KL$ to simplify notation). With increasing time, the widths of the triangle-shaped peaks decrease but their maxima at $x = \pm 5$ 
become higher while unchanged at $x = 0$. At sufficiently long times $t \gtrsim 2 \cdot 10^4$ 
the modification of the density profile slows down to a level suggesting that the
system approaches a non-trivial stationary state in which $\partial n(x,t)/ \partial t = 0$. 
In other words, the shape of the density profile is transformed to such a form at which 
any coalescence processes become impossible in view of the specific form of the coalescence kernel 
$b(x) = C_{1,1,8} (x)$. The latter accepts nonzero values only in the interval 
$s' - \sigma = 7 < |x| < 9 = s' + \sigma$, so that the absence of coalescence at a given 
spatial configuration means that there exists no pair of particles with interparticle separations 
$|x|$ lying in the interval $[7, 9]$. Allowing particles to jump changes the situation radically, 
as is demonstrated in Fig. \ref{fig:3}(b) for Gaussian jump kernel $G_{0.2,1}(x)$. Even for relatively 
small jump intensity $\mu_a = 0.2$ and range $\sigma_a = 1$, the density quickly decreases to zero 
for each $x$, after initial period of time, when new peaks are formed. It is interesting to remark 
that the monotonic decrease of main maxima in $x \pm 5$ is accompanied by nonmonotonic change of the 
magnitude of the newly formed peaks in $x = 0$ (and $\pm KL$). This magnitude first increases at 
$0 < t \lesssim 4$ achieving a maximum at $t \sim 4$, and then decreases at $t \gtrsim 4$.

\begin{figure}
	\centering
	\includegraphics{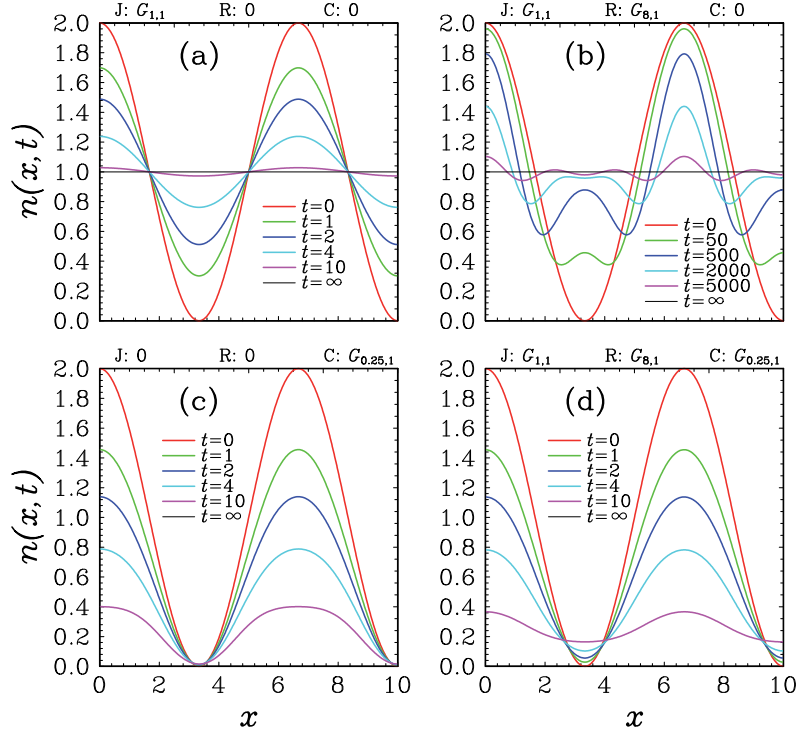}
	\caption{
		Density profile $n(x,t)$ starting from trigonometric initial condition $n(x, 0) = T_{1,1,3} (x)$. 
		The jump, repulsion, and coalescence kernels are Gaussians $G_{1,1}(x)$, $G_{8,1}(x)$, 
		and $G_{0.25,1} (x)$, respectively. Other notations	are the same as for Fig. \ref{fig:1}.
	}
	\label{fig:4}
\end{figure}

{\it Trigonometric initial density profiles.} --- Another interesting case is the initial density in
the form of a trigonometric function,
\begin{equation}\label{eq:Trigonometric}
	n(x,0) = T_{n_0, \mu_0, k} (x) = n_0 \Big( 1 + \mu_0 \cos (2 \pi k x / L) \Big),
\end{equation}
where $0 < \mu_0 < 1$ is the coefficient of the modulation and $k \geq 1$ defines the period 
$L/k$ in coordinate. Then PBC should be used to reproduce the infinite system. The obtained
solutions $n(x,t)$ for $n(x, 0) = T_{1,1,3} (x)$ with $n_0 = 1$, $\mu_0 = 1$, $k = 3$, and 
$L = 20$ are plotted in Fig. \ref{fig:4} when jump, repulsion, and coalescence kernels are 
Gaussians $a(x) = G_{1,1} (x)$, $\varphi(x) = G_{8,1} (x)$, and $b(x) = G_{0.5,1} (x)$, respectively. 
For convenience, we presented $n(x,t)$ only on the right-hand side of coordinate space at $x \geq 0$.

From Fig. \ref{fig:4}(a) we see that pure free jumps do not change the form of density profile
which remains to be of the trigonometric shape with the same frequency $k$ and periodicity $L$.
In particular, the density continues to oscillate around the same level $n_0 (1 + \mu_0 )/2 = 1$ 
for any t. However, the amplitude of these oscillations decreases with increasing $t$, 
so that the particle density profile seems to become flat at long times, 
$\lim t \rightarrow \infty n(x,t) = 1$. Again, the total number $\int_{\Omega} n(x,t)dx$ of particles 
in $\Omega \in [-L/2, L/2]$ does not change in time. The repulsion makes the simple trigonometric 
(cosine) form much more complicated with the existence of additional maxima and minima inside 
$\Omega$, see Fig. \ref{fig:4}(b). Moreover, the homogeneity is being achieved here much slower than 
in the case of free jumps (compare density at $t \gtrsim 5000$ versus for $t \gtrsim 10$ in 
Fig \ref{fig:4}(a)). This can be explained by the strong intensity ($\mu_\varphi = 8$) of
repulsion potential $\varphi(x) = G_{8,1} (x)$.

For pure coalescence in Fig. \ref{fig:4}(c), the distribution $n(x,t)$ is no longer of a trigonometric
form at $t > 0$ contrary to pure free jumps and like as in the case of repulsion. In addition, here
the density seems to decrease to zero at $t \rightarrow \infty$. The inclusion of repulsive jumps changes
the behaviour of $n(x,t)$ near its minima in a characteristic way. Indeed, in Fig. \ref{fig:4}(c) these
minima remain to be equal to zero, while in Fig. \ref{fig:4}(d) they have a tendency to increase to
positive values at some $t$. On the other hand, the speed of decrease of maxima in $n(x,t)$
with increasing time practically does not change at $t \lesssim 10$ despite their strength. Moreover,
the shape of the density profile at long $t \gtrsim 10$ is modified as well, making it more flat with
smaller oscillations. As a result, spatial homogeneity is obtained here faster.

{\it Step initial density profiles.} --- A special case presents initial condition in the form of
the Heaviside step function
\begin{equation}\label{eq:Heaviside}
	n(x,0) = H_{n_0}(x) = \left\{
	\begin{array}{rcl}
		n_0, & & x \leq 0\\
		0, & & x > 0
	\end{array}
	\right. .
\end{equation}
Here the system is initially ($t = 0$) considered on the finite interval $[-L/2, L/2]$ with no
PC. Further ($t > 0$) its size is gradually increased to the infinity on the unbounded interval
$x \in ]-\infty, \infty[$ at $t \rightarrow \infty$ according to the automatically adjusted 
system-size (AASS) approach.
Then, a modified BC should be applied by combination of DBC and ABC together with
AASS when analyzing the densities on the boundaries $\pm L/2$. The DBC are used from the
right, where $\lim_{x \rightarrow \infty} n(x,t) = 0$ for all $t$. From the left, where 
$\lim_{x \rightarrow - \infty} n(x,t) = n(t) \neq 0$ with $n(0) = n_0$, we must employ ABC. 
This means that we should exploit the DBC rule for $n_{i+j}$ given by the second equality of 
Eq. \eqref{eq:DBC}, i.e., replace $n_{i+j}$ by $0$ if $i + j > N$, and
the ABC rule for $n_{i-j}$ given by the first equality of Eq. \eqref{eq:ABC}, i.e., replace 
$n_{i-j}$ by $n_1$ if $i - j < N$. When nonzero values of $n(x,t)$ approach the right boundary 
at $x = L/2$, i.e., when $n(L/2,t) > \varepsilon \max_x n(x,t)$, the system size $L$ is enlarged 
(in two times) and the simulations are continued.

By monitoring from the left the difference between the actual values of $n(x,t)$ at $x =-L/2$ and their 
homogeneous counterpart $n^h_{RK} (t)$ [obtained by solving numerically the kinetic equation for the 
spatially homogeneous initial condition $n(x, 0) = n_0$ in parallel to the spatially inhomogeneous case 
using the same RK method] we can estimate the influence of the finiteness of $L$. If this difference 
exceeds the predefined level $\varepsilon \max_x n(x,t)$
we should enlarge $L$ according to the automatically adjustable system-size approach. The
asymptotic value $n_1 (t)$ will differ from the exact solution $n^h (t)$ even at very large $L$ because 
of the approximate character of time integration. In the limit of sufficiently small
step sizes $\Delta t$ when $\lim_{\Delta t \rightarrow 0, L \rightarrow \infty} n_1 (t) = n^h (t)$ 
we come to the limiting ABC-DBC scheme:
$\lim_{x \rightarrow - \infty} n(x,t) = n^h (t)$ and $\lim_{x \rightarrow \infty} n(x,t) = 0$. 
Calculating the difference between the actual values of $n(x,t)$ at $x = -L/2$ and their exact 
counterpart $n^h (t)$ at $x \rightarrow - \infty$ we can
estimate the influence of two effects in one fashion, namely, those caused by the finiteness
of $L$ (should be significantly large) and $\Delta t$ (should be significantly small). 
In such a way, both ABC and DBC deviations are analyzed when deciding on the size enlargement within
the ABC-DBC scheme. This completes the automatically adjustable system-size approach.

\begin{figure}
	\centering
	\includegraphics{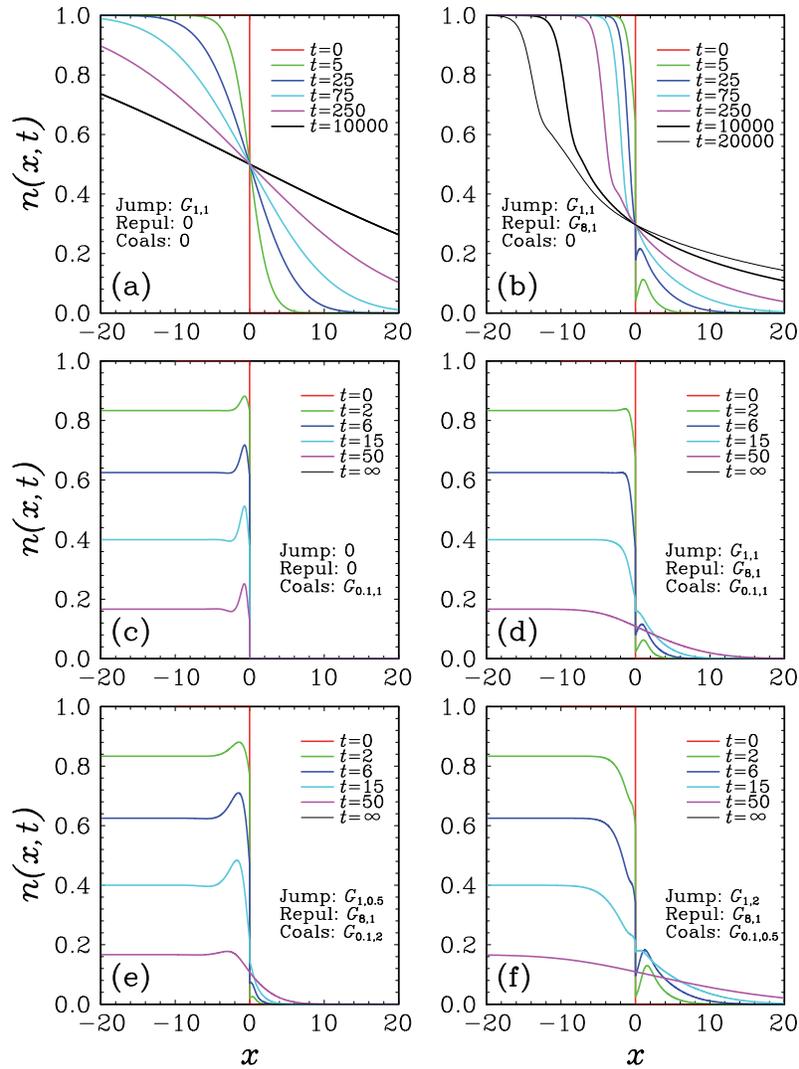}
	\caption{
		Time evolution of density profile for initial condition $H_1 (x)$. The jump, repulsion, 
		and coalescence	kernels are Gaussians with different intensities and ranges 
		(see the legends inside). Initially ($t = 0$) the system is considered on the finite interval 
		$[-10, 10]$ with no periodic conditions and further ($t > 0$) its size gradually
		increases to the infinity on $x \in ]-\infty, \infty[$ at $t \rightarrow \infty$ 
		according to the automatically adjusted approach. Other	notations are the same as for 
		Fig. \ref{fig:1}.
	}
	\label{fig:5}
\end{figure}

Time evolution of $n(x,t)$ for $n(x, 0) = H_1 (x)$ is presented in Fig. \ref{fig:5} using Gaussian jump
$G_{\mu_a,\sigma_a}$, repulsion $G_{\mu_\varphi, \sigma_\varphi}$, and coalescence $G_{\mu_a, \sigma_b}$ 
kernels with different intensities of $\mu_a = 1$, $\mu_\varphi = 8$, and $\mu_b = 0.1$, respectively, 
as well as with equal [$\sigma_{a,\varphi,b} = 1$, parts (a)--(d)] or different 
[$\sigma_a = 0.5$, $\sigma_\varphi = 1$, $\sigma_b = 2$, part (e), and 
$\sigma_a = 2$, $\sigma_\varphi = 1$, $\sigma_b = 0.5$, part (f)] ranges.
As can be seen from Fig. \ref{fig:5}(a) for pure free jumps, the discontinuous step function 
$n(x, 0) = H_1 (x)$ transforms into a continuous S-shaped curve at finite times $t > 0$. All the curves
intersect each other in the same point $(0, 1/2)$, where $1/2$ is the arithmetic mean of two
initial values ($1$ to the left and $0$ to the right). The slope of these curves becomes smaller
with increasing time, so that the system tends to the mid-value everywhere. Moreover, the
decrease of the amount of particles for $x < 0$ is equal to the corresponding increase for
$x > 0$, that can be written as $\int_{-\infty}^{0}[1 - n(x,t)]dx = \int_0^\infty n(x,t) dx$. 
The same statement concerns repulsive jumps, but here the curves are intersected in a point 
which lies below $(0, 1/2)$.
In addition, at short times $t \lesssim 25$ the density profile $n(x,t)$ remains to be discontinuous in
$x = 0$, and the shape of the curves is more complicated, including the existence of maximum
in $x \sim 1$. The latter disappears at $t \gtrsim 25$, and $n(x,t)$ becomes more and more flat with the
flow of time. However, the flattening process is not as fast as in the case of free jumps.

For pure coalescence we can observe in Fig. \ref{fig:5}(c) that the initial step function 
$n(x, 0) = H_1 (x)$ changes to a step-like dependence $n(x,t)$ at $t > 0$ containing one small peak in 
$x \sim -1$. The density in $]-\infty, 0[$ decreases from $1$ at $t = 0$ to zero at larger $t > 0$. 
No particles appear at all in $]0, \infty[$ for any $t$. Moreover, the initial discontinuity does 
not vanish even for relatively long times. 
The inclusion of repulsive jumps, see Fig. \ref{fig:5}(d), prevents the appearance of peaks at 
$x \sim -1$, while at $x > 0$ the profile $n(x,t)$ is more flat when compared to
that in Fig. \ref{fig:5}(b). Here $n(x,t)$ decreases to zero at sufficiently large times 
$(t \gtrsim 100)$ as in the case of pure coalescence. As the range of jumps becomes shorter, 
$\sigma_a = 0.5$, and the range of coalescence increases, $\sigma_b = 2$, see Fig. \ref{fig:5}(e), 
the peaks in $n(x,t)$ appear again at $x \sim -\sigma_b = -2$ (i.e. they shift to the left with 
respect to those for $\sigma_b = 1$). Moreover, they become more pronounced with larger amplitudes. 
The right lying peaks at $x \sim 1$ shifts their positions to $x \sim 0.5$ and decrease their 
amplitudes. This is contrary to the case when the jump range increases, $\sigma_a = 2$, 
and the coalescence range decreases, $\sigma_b = 0.5$, see Fig. \ref{fig:5}(f).
Then the peaks on the left do not appear, while the maximum on the right becomes more
visible. In addition, the density approaches zero more slowly here. Note that for the inverse 
initial unit step function $n(x, 0) = H_1 (-x)$, the results $n(x,t)$ presented in Fig. \ref{fig:5} for
$n(x, 0) = H_1 (x)$ should also be inversed by $n(-x,t)$ to obtain solutions without resolving
the kinetic equation. This statement is quite general and remains in force not only for step
functions, but for any other asymmetric initial conditions. Namely, when picking the initial condition 
$n^* (x, 0) = n(-x, 0)$, we automatically obtain $n^* (x,t) = n(-x,t)$, where $n(x,t)$
is the solution for $n(x, 0)$. This follows from the structure of kinetic equation 
\eqref{eq:kinetic} and the symmetricity of kernel functions.

\begin{figure}
	\centering
	\includegraphics{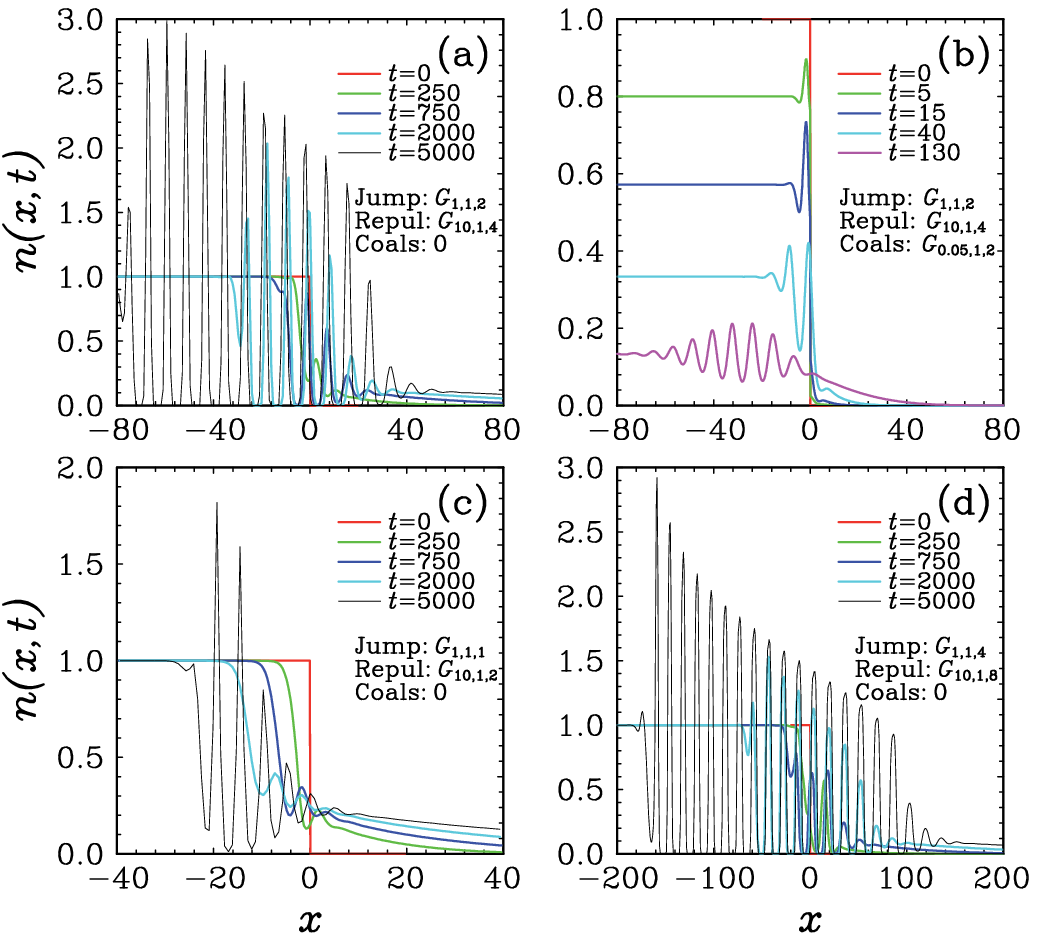}
	\caption{
		Dynamics of the system starting from unit step function $H_1(x)$ in the cases of pure 
		repulsive jumps
		with different kernel shifts: (a) $s = 2$, $s' = 4$; (c) $s = 1$, $s' = 2$; and 
		(d) $s = 4$, $s' = 8$; as well as of (b)
		repulsive jumps and coalescence with $s = 2$ and $s' = 4$. Initially ($t = 0$) the system 
		is considered on the finite
		interval $[-20, 20]$ with no periodic conditions and further ($t > 0$) its size gradually 
		increases to the infinity on
		$x \in ]-\infty, \infty[$ at $t \rightarrow \infty$ according to the automatically adjusted 
		approach. The interactions are described by the
		shifted pair Gaussian jump $G_{1,1,s}$, repulsion $G_{10,1,s'}$ and coalescence 
		$G_{0.05,1,2}$ kernels (see the legends inside).
	}
	\label{fig:6}
\end{figure}

Consider, finally, the dynamics of the initial step distribution $n(x, 0) = H_1 (x)$ in the presence 
of pairs of shifted Gaussian kernels. First, we used shifting parameter $s = 2$ for moderate 
($\mu_a = 1$) jump kernel $a(x) = G_{1,1,2} (x)$ and weak ($\mu_b = 0.05$) coalescence interaction
$b(x) = G_{0.05,1,2} (x)$, while $s' = 4$ for strong ($\mu_\varphi = 10$) repulsion potential 
$\varphi(x) = G_{10,1,4} (x)$ (the ranges were $\sigma_a = \sigma_b = \sigma_\varphi = 1$ for all 
the kernels). The corresponding results are presented in Fig. \ref{fig:6}. From part (a) of this figure 
we see that for pure repulsive jumps such a choice leads at $t > 0$ to the emergence of 
self-propagating spatial inhomogeneity in the from of persistent oscillations with thin peaks of 
width of order of $s = 2$ at level of $n(x) = 1$ 
and distance between them of order of $2s' = 8$. For $n(x, 0) = H_1 (x)$ the inhomogeneous front
propagates to the left with increasing amplitude at not too large $x$ and to the right with 
dumping oscillations (for $n(x, 0) = H_1 (-x)$ the pattern is inverse). The inclusion of coalescence
even with a slight intensity of $\mu_b = 0.05$ drastically changes the situation, see Fig. \ref{fig:6}(b).
Here, the oscillations are not so strong, especially at $x > 0$, and they quickly disappear with
increasing time.

Further, we decreased and increased the shifting parameters in two times to values $s = 1$
and $s' = 2$ as well as to $s = 4$ and $s' = 8$. The results for these two cases are shown in parts (c)
and (d) of Fig. \ref{fig:6}, respectively, at the absence of coalescence. As can be seen from 
Fig. \ref{fig:6}(c), the decrease of $s$ and $s'$ suppresses the processes of density propagation 
by lowering its speed and oscillation amplitude. Moreover, the distance between peaks in $n(x,t)$ 
and their width also decrease correspondingly to $2s' = 4$ and $s = 1$. This is in a contrast to 
the opposite case when the kernel shifts increase. Such an increase leads to stimulation of 
density propagation with high amplitude of the oscillations. Then the distance between peaks and their 
width increase to $2s' = 16$ and $s = 4$, respectively. The propagation speed grows proportionally as
well.

{\it Precision of the simulations.} — Accuracy of our numerical simulations was measured
in terms of relative absolute deviations of density values $n_i = n(x_i, t)$ in grid points $x_i$ from
``exact'' data $\breve{n}_i$ at time $t$ using the relation
\begin{equation}\label{eq:Accuracy}
	\Theta(h, \Delta t) = \frac{\sum_{i=1}^{\breve{N}} |n_i - \breve{n}_i|}{\sum_{i=1}^{\breve{N}} \breve{n}_i}
	\times 100 [\%].
\end{equation}
The total number $\breve{N} \leq N$ of grid points involved into summation \eqref{eq:Accuracy} 
depends on the spatial region considered. The ``exact'' (or rather reference) values $\breve{n}_i$ 
were obtained at high enough
space and time resolutions with a tiny mesh of $h = 0.005$ and a tiny time step of $\Delta t = 0.01$
in order to be entitled to ignore the numerical uncertainties. The spatial and time integrations 
were performed with the help of the composite Simpson rule and RK4 algorithm,
respectively. The numerical error analysis was done by carrying out a series of simulations
at different $h = 0.01, 0.02, 0.04, 0.1, 0.2$ and $\Delta t = 0.1, 0.2, 0.4, 1.0, 2.0, 2.5$. 
The numerical deviations were then estimated by Eq. \eqref{eq:Accuracy} to plot 
$\Theta (h, \Delta t)$ in a wide range of varying $h$ and $\Delta t$. 
For the purpose of comparison, the simulation runs with the composite trapezoidal rule
and RK2 algorithm were performed as well. Error results presented below are related to one
of the situations considered in Figs. \ref{fig:1}--\ref{fig:6} when choosing initial condition 
$n(x, 0)$, namely, to the case of Fig. \ref{fig:5}(d) at $t = 50$. There, the homogeneous 
and inhomogeneous intervals were $x \in [-80, -20]$ and $x \in [-20, 20]$, respectively. Similar 
results were observed for all the rest choices of $n(x, 0)$.

\begin{figure}
	\centering
	\includegraphics{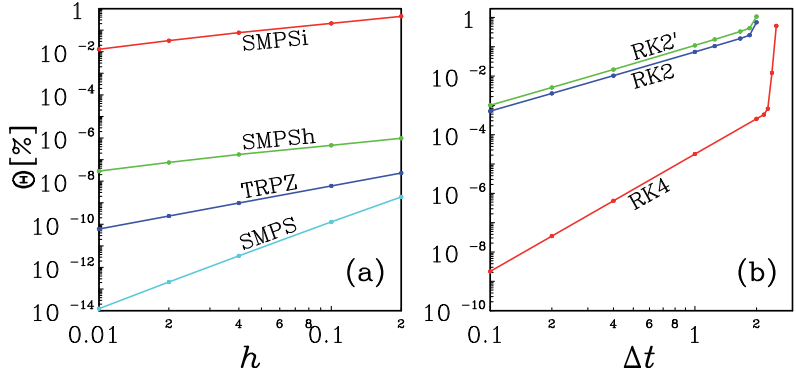}
	\caption{
		(a) Uncertainties $\Theta$ of the simulations related to single composite Simpson (SMPS) 
		or trapezoidal (TRPZ) spatial integrations and the overall spatial discretization obtained 
		in homogeneous (SMPSh) and inhomogeneous (SMPSi) regions at different mesh $h$; 
		(b) Dependencies of $\Theta$ on the size of the time step $\Delta t$
		for the RK2, RK2', and RK4 time integrations.
	}
	\label{fig:7}
\end{figure}

The numerical uncertainties $\Theta$ cased by spatial discretization and spatial integration are
shown in Fig. \ref{fig:7}(a) as functions of the length $h$ of space mesh (at a fixed RK4 time step
of $\Delta t = 0.01$). Note that the log-log presentation was employed to show the behaviour of
$\Theta (h, \Delta t)$ at small $h$ and $\Delta t$ in more detail. Here we distinguish three 
kinds of errors. The first is related to uncertainties due to single spatial integrations. 
The corresponding results
for them obtained by the composite Simpson and composite trapezoidal rules are marked
by SMPS and TRPZ, respectively. We see that for sufficiently small values of $h$, the relative 
SMPS or TRPZ errors $\Theta$ are proportional to $h^4$ or $h^2$. Indeed, the log-log plots are
lines with slopes 2 or 4, i.e., $\log \Theta = c_2 + 2 \log h$ or $\log \Theta = c_4 + 4 \log h$, 
where $c_2$ and $c_4$ denote some constants. 
This confirms the fact \cite{NumericalMethods} that the composite Simpson and trapezoidal rules 
are accurate to the fourth $\mathcal{O}(h^4)$ and second $\mathcal{O}(h^2)$ orders in $h$, respectively. 
The accuracy is lowered significantly when considering the overall uncertainties caused by spatial discretization of the kinetic equation. Such uncertainties observed in the homogeneous
and inhomogeneous interval with the Simpson integration are presented by curves labeled
correspondingly as SMPSh and SMPSi. The TRPZh and TRPZi data of the trapezoidal integration 
appear to be approximately the same (see below the reason) and are not shown in the figure.

Maximal uncertainties are achieved in the inhomogeneous region because the dependence of 
$n(x,t)$ on $x$ makes the functions under spatial integrals more sharp, decreasing
the precision of the discretization. The local uncertainties of such a discretization are of
order of $\mathcal{O}(h^3)$. Taking into account that the total number of numerical operations 
is proportional to $N^2$, the global errors appear as an accumulation of local ones leading to the
overall uncertainties $\Theta$ of order of $N^2 \mathcal{O}(h^3) \sim \mathcal{O}(h)$, 
since $h = L/N$. In other words, the SMPSh and SMPSi deviations behave as 
$\log \Theta = c_1 + \log h$. Despite the linear dependence $\mathcal{O}(h)$, even the SMPSi 
errors are relatively small and do not exceed about $0.01$ -- $0.03\%$ at
$h = 0.01$ -- $0.02$ (see Fig. \ref{fig:7}). They are, however, much larger than those of the 
single spatial integration. Therefore, the errors caused by spatial discretization of the kinetic 
equation significantly prevail over the spatial integration uncertainties (the both use the same mesh $h$).
For this reason it is not so important what the method (Simpson or trapezoidal) is applied to
spatial integration. This explains why the SMPSh/SMPSi and TRPZh/TRPZi data are very
close to each other.

Dependencies of $\Theta$ on the size of the time step $\Delta t$, obtained in the simulations with
the RK2, RK2', and RK4 time integrations (at a fixed spatial mesh of $h = 0.01$ with the
involved inhomogeneous interval $x \in [-20, 20])$ are depicted in Fig. \ref{fig:7}(b). It can be seen
clearly that for sufficiently small values of $\Delta t$, the relative deviations $\Theta$ are proportional to
$\Delta t^2$ or $\Delta t^4$ for the algorithms of the second (RK2 and RK2') or fourth (RK4) orders, respectively. Here, the log-log plots are lines with slopes 2 or 4, i.e., 
$\log \Theta (\Delta t) = c_2 + 2 \log \Delta t$ or $\log \Theta (\Delta t) = c_4 + 4 \log \Delta t$ 
[similarly as for dependence of $\Theta$ on $h$ for the single Simpson and
trapezoidal spatial integrations, see Fig. \ref{fig:7}(a)]. We mention that the local (single-step) errors
of the RK4 algorithm are of order of $\mathcal{O}(\Delta t^5)$, see Eq. \eqref{eq:RK4}. 
So that the numerical integration over long times $t \gg \Delta t$ leads to the global errors of order 
of $t / \Delta t \mathcal{O} (\Delta t^5) = \mathcal{O} (\Delta t^4)$ as this
requires repeating of the single-step integration by $t / \Delta t$ times. Analogously, the RK2 and
RK2' local uncertainties $\mathcal{O}(\Delta t^3)$ modifies to the global $\mathcal{O}(\Delta t^2)$-errors. 
Since $\Delta t^4$ decreases with decreasing $\Delta t$ much faster than $\Delta t^2$, the 
fourth-order RK4 algorithm produces much
more precise solutions than the second-order integrators RK2 and RK2'. Thus, the former
can be recommended when a very high accuracy is required. For instance, the RK4 algorithm is precise up to a negligible level of order of $\Theta \sim 10^{-9}\%$ at $\Delta t \sim 0.1$. On the other
hand, the RK2 and RK2' integrators can provide here only an accuracy of $10^{-3}$ [RK2 is
slightly better than RK2', see Fig. \ref{fig:7}(b)]. Remember, however, that the overall errors include
also the spatial discretization uncertainties which are of order of $10^{-2}\%$ at $h = 0.01$, see
Fig. \ref{fig:7}(a). For this reason, the RK4 algorithm can be applied with much larger time steps of
$\Delta t \sim 1$ -- $2$, where the time-discretization uncertainties did not exceed $\sim 10^{-3}$ 
and are comparable with the space-discretization errors. This significantly accelerates the simulations
since the computational time is inverse proportional to $\Delta t$. In particular, the RK4 speedup
at $t = 1$ (when $\Theta \sim 10^{-5} \%$) with respect to the RK2 and RK2' schemes at $\Delta t = 0.1$ 
(where a worse precision of $\Theta \sim 10^{-3}\%$ is observed despite the smaller time step) is of 
order of $1/0.1/2=10/2=5$ (the factor ``2'' takes into account the fact that RK4 requires a twice larger
number of operations per step than those of RK2 or RK2' ). With further increasing $\Delta t$, all
the integrators become unstable and cannot be used at $\Delta t > 2$, where $\Theta > 1\%$.

\section{Conclusion}
In this paper we have derived an efficient algorithm to obtain numerical solutions for the
time-differential kinetic equation which approximates nonlocal stochastic evolution of coalescing and repulsively jumping particles in the continuous space $\mathbb{R}^d$. The equation is very
difficult from the analytical point of view due to the presence of complicated spatial integrals with nonlinear and nonlocal terms and therefore requires numerical analysis. The
proposed algorithm is based on a set of techniques including time-space discretization, periodic, Dirichlet, and asymptotic boundary conditions, composite Simpson and trapezoidal
rules, Runge-Kutta schemes, and automatically adjustable system-size approaches. This has
allowed as to carry out simulations of dynamics for one-dimensional systems with various
initial inhomogeneous densities and different forms of the coalescence, jump and repulsion
kernels, giving a comprehensive study of the jump-coalescence dynamics. A numerical error
analysis of the obtained results has also been performed.

For some specific choices of the model parameters and initial densities, a nontrivial
dynamics has been revealed. In the case of pure free jumps the system always tends to a
nonzero homogeneous density for any initial conditions. In contrast, the presence of strong
repulsion potentials can result in the appearance of persistent wave-like density propagations
when the repulsion and jump kernels are chosen in the form of a sum of shifted single
Gaussian functions. The shifting parameters define the shape and spatial period of density
structures. The introduction of even relatively weak coalescence prevents them to persist.
In the case of pure coalescence, the population eventually goes extinct, except for a special
choice of the initial density profile and coalescence kernel. For instance, if the particles are
initially located exclusively on an archipelago of islands and the minimal distance between
them is equal to the shifting parameter of the coalescence kernel consisting of two simple
rectangle functions with the ranges which do not exceed the sizes of the islands, then a
stationary state with inhomogeneous particle distributions can arise. The inclusion of any
jumps even with extremely small intensity radically changes the situation, leading to the
collapse of the system.

The proposed algorithm is implemented in a Fortran program code which can be received by request and 
exploit by any researcher free of charge. The code can be adapted
to more complex multicomponent models of jumping and coalescence particles of different
types. The numerical simulations can be extended to systems of higher dimensions. The
Poisson approximation used for the derivation of the kinetic equation can be improved by
advancing to the Kirkwood \cite{BeyondMeanfield} or Fisher-Kopeliovich \cite{BeyondKirkwood} 
ansatz like for birth-death models. Mass and size of particles can also be taken into account. 
These and other topics as well
as possible applications of the obtained results to real populations will be the subject of our
further investigations.

\section*{Appendix}
\renewcommand{\theequation}{A\arabic{equation}}
\setcounter{equation}{0}
In the absence of coalescence [$b(x) = 0$] we can apply the convolution theorem to avoid the
direct spatial integration. Indeed, then the kinetic equation \eqref{eq:kineticDiscrete} can 
be rewritten in the form
\begin{equation}\label{eq:kineticConv}
	\frac{\partial n(x, t)}{\partial t} = - n(x, t) \int a(x - y) g(y, t) dy
	+ g(x, t) \int a(x - y) n(y, t) dy,
\end{equation}
where
\begin{equation}\label{eq:gConv}
	g(x, t) = \exp \Big( - \int \varphi (x - u) n(u, t) du \Big) .
\end{equation}
To simplify notation consider the case $d = 1$ and introduce the discrete direct and inverse
one-dimensional spatial Fourier transforms \cite{NumericalRecipes}
\begin{equation}\label{eq:Fourier}
	\tilde{f}_k = \sum_{i=0}^{N-1} f_i \exp (-2 \pi \imath k i / N) \equiv \{f_i\}_k, \quad
	f_i = \frac{1}{N} \sum_{k=0}^{N-1} \tilde{f}_k \exp (2 \pi \imath k i / N) \equiv \{f_k\}_i^{-1}
\end{equation}
where $f_i = f(x_i)$ and $x_i = -L/2 + (i - 1/2)h \in ]-L/2, L/2[$ with even $N$ and $h = L/N$. 
Now, applying the convolution theorem, the discrete counterpart of Eq. \eqref{eq:kineticConv} 
in view of Eq. \eqref{eq:gConv} can be presented as
\begin{equation}\label{eq:derivativeConv}
	\dot{n}_i = \frac{dn_i}{dt} = -h n_i \{\tilde{a}_k \tilde{g}_k\}_i^{-1} + h g_i \{\tilde{a}_k \tilde{n}_k\}_i^{-1}, \quad g_i = \exp (-h \{ \tilde{\varphi}_k \tilde{n}_k \}_i^{-1})
\end{equation}
where $n_i = n(x_i ,t)$ with $\tilde{n}_k$, $\tilde{a}_k$, and $\tilde{\varphi}_k$ being the discrete 
Fourier transforms of $n_i, a_i = a(x_i)$, and $\varphi_i = \varphi(x_i)$ according to Eq. 
\eqref{eq:Fourier} at $f_i \equiv n_i, a_i, g_i$ and 
$\tilde{f}_k \equiv \tilde{n}_k, \tilde{a}_k, \tilde{\varphi}_k$.

Because the kernel function $a(x)$ and repulsive potential $\varphi(x)$ do not depend on time,
their Fourier components $\tilde{a}_k$ and $\tilde{\varphi}_k$ can be calculated once at the very 
beginning. Further, having the current values $n_i$ for $i = 0, 1, \ldots, N - 1$ we compute their 
Fourier duplicates $\tilde{n}_k$ for $k = 0, 1, \ldots, N - 1$. On the basis of the known products 
$\tilde{a}_k n_k$ and $\tilde{\varphi}_k \tilde{n}_k$ we evaluate their coordinate counterparts 
$\{\tilde{a}_k \tilde{n}_k \}_i^{-1}$ and $\{\varphi_k N_k\}_i^{-1}$ by exploiting the inverse 
Fourier transform.
Constructing $g_i = \exp(-\{\tilde{\varphi}_k \tilde{n}_k \}_i^{-1})$ we calculate $\tilde{g}_k$ 
by means of the direct Fourier transform and use again the inverse transform to obtain 
$\{\tilde{a}_k \tilde{g}_k \}_i^{-1}$. In such a way we have all
necessary components to form the time derivatives $\dot{n}_i$ according to Eq. 
\eqref{eq:derivativeConv}. Despite this requires several direct and inverse transformations, 
the order of total number of operations
can be reduced from $N^2$ (when the spatial integration is performed directly, see Sect. 
\ref{sec:Algorithm}) to $N \ln N$ when exploiting the fast discrete direct and inverse 
Fourier transforms. This is very important feature because $N \gg 1$ is large (in practice 
$N \sim 10^2$ -- $10^3$ at least). However, such an approach can work not so well for 
discontinuous functions. Moreover, it assumes that all involved functions are periodic on 
the interval $[-L/2, L/2]$.

The generalization to any higher dimensionality $d > 1$ can also be done by replacing $i$
and $k$ on the $d$-dimensional vectors $\mathbf{i} = (i_1, i_2, \ldots, i_d)$ and 
$\mathbf{k} = (k_1, k_2, \ldots, k_d)$ as well as their multiplication $ki$ on the scalar 
product $\mathbf{k} \cdot \mathbf{i}$ with summation in Eq. \eqref{eq:Fourier} from $0$ to $N - 1$ for
each dimension, where $f_\mathbf{i} = f(x_{i_1}, x_{i_2}, \ldots, x_{i_d})$. As a result, the total 
number of operations increases from $N^2$ to $N^{2d}$ or from $N \ln N$ to $N^d \ln N^d$ when applying 
the product of $d$ one-dimensional usual or fast Fourier transform in Cartesian coordinates. 
If the functions possess radial symmetry, we can consider the Fourier transform in polar or 
spherical coordinates at $d = 2$ or $3$, for instance. This allows to reduce $d$-dimensional 
integrations to discrete Hankel \cite{Convolution2D} or Fourier-Bessel \cite{Convolution3D} transforms 
in one (radial) dimension by integrating analytically out all angle dependencies. The chief advantage 
of this approach is the fact that then the total number of operations becomes independent of $d$ 
and remains the same (of order of $N \ln N$) for any dimensionality \cite{Nikitin}, just as at $d = 1$.

In the spatially homogeneous case, i.e. when $n(x,t) = n(t)$ does not change on $x$, we can
carry out the spatial integrations in Eqs. \eqref{eq:kineticDiscrete} and \eqref{eq:repulsionDiscrete} 
over $y$ explicitly. This significantly simplifies the kinetic equations to the form 
$dn/dt = -\mu_b n^2$ which gives the analytical solution $n(t) = n(0)/[1 + \mu_b n(0)t]$. 
It tends to zero with $t \rightarrow \infty$ for any initial density $n(0)$ and $\mu_b \neq 0$.
In this case the jump contribution completely disappears. At $\mu_b = 0$ this leads to the 
time-independent solution $n(t) = n(0)$. It disappears also in the spatially inhomogeneous case
for systems with finite number $\mathcal{N} = \int n(x)dx$ of particles when calculating 
$d \mathcal{N} /dt$ (then during integration in the rhs of Eq. \eqref{eq:kinetic} over $x$, 
the first two terms are mutually cancelled).


%
%


\begin{thebibliography}{}
%
%
\bibitem{Arratia} 
Arratia, R.A.: Coalescing Brownian motion on the line. PhD thesis, University of Wisconsin, 
Madison, ProQuest LLC, Ann Arbor, MI (1979)
\bibitem{Toth}
T\'oth, B., Werner W.: The true self-repelling motion. Probab. Theory Relat. Fields 111, 
375--452 (1998)
\bibitem{LeJan}
Le Jan, Y., Raimond, O.: Flows, coalescence and noise. Ann. Probab. 32, 1247--1315 (2004)
\bibitem{Konarovskii}
Konarovskii, V.V.: On an infinite system of diffusing particles with coalescing. Teor. Veroyatn. 
Primen. 55, 157--167 (2010)
\bibitem{Berestycki}
Berestycki, N., Garban, Ch., Sen, A.: Coalescing Brownian flows: a new approach. Ann. Probab. 43,
3177--3215 (2015)
\bibitem{Renesse}
Konarovskii, V.V., von Renesse, M.: Modified massive Arratia flow and Wasserstein diffusion. Comm.
Pure Appl. Math. 72, 764--800 (2019)
\bibitem{RCRJ1}
Pilorz, K.: A kinetic equation for repulsive coalescing random jumps in continuum. Ann. Univ. 
Mariae Curie-Sk{\l}odowska Sect. A 70, 47--74 (2016)
\bibitem{RCRJ2}
Kozitsky, Yu., Pilorz, K.: Random jumps and coalescence in the continuum: evolution of states of an
infinite system. ArXiv:1807.07310 (2018)
\bibitem{Baranska}
Bara\'nska, J., Kozitsky, Yu.: The global evolution of states of a continuum Kawasaki model with 
repulsion. IMA J. Appl. Math. 83, 412--435 (2018)
\bibitem{Berns}
Berns, C., Kondratiev, Y., Kozitsky, Y., Kutoviy, O.: Kawasaki dynamics in continuum: micro- and 
mesoscopic descriptions. J. Dynam. Differential Equations 25, 1027--1056 (2013)
\bibitem{Delius}
Capit\'an, J.A., Delius, G.W.: Scale-invariant model of marine population dynamics. Phys. 
Rev. E 81, 061901 (2010)
\bibitem{Blanchard}
Law, R., Plank, M.J., James A., Blanchard J.L.: Size-spectra dynamics from stochastic predation and
growth of individuals. Ecology 90, 802--811 (2009)
\bibitem{Weisse}
Weisse, T.: Dynamics of autotrophic picoplankton in marine and freshwater ecosystems. In: Jones, 
J.G.(ed.) Advances in Microbial Ecology, vol. 13, pp. 327--370. Plenum Press, New York (1993)
\bibitem{RudnickiWieczorek1}
Rudnicki, R., Wieczorek, R.: Fragmentation-coagulation models of phytoplankton. Bull. Pol. Acad. 
Sci. Math. 54, 175--191 (2006)
\bibitem{RudnickiWieczorek2}
Rudnicki, R., Wieczorek, R., Phytoplankton dynamics: from the behaviour of cells to a transport 
equation. Math. Model. Nat. Phenom. 1, 83--100 (2006)
\bibitem{Ambrose}
Ambrose J., Livitz M., Wessels D., Kuhl S., Lusche D.F., Scherer A., Voss E., Soll D.R.: Mediated
coalescence: a possible mechanism for tumor cellular heterogeneity. Am. J. Cancer Res. 5, 3485--3504
(2015)
\bibitem{Wessels}
Wessels D., Lusche D.F., Voss E., Kuhl S., Buchele E.C., Klemme M.R, Russell K.B., Ambrose J., Soll
B.A., Bossler A., Milhem M., Goldman C., Soll D.R.: Melanoma cells undergo aggressive coalescence
in a 3D Matrigel model that is repressed by anti-CD44. PLoS One. 12, e0173400 (2017)
\bibitem{RCRJ3}
Kozitsky, Yu., Omelyan I., Pilorz K.: Jumps and coalescence in the continuum: a numerical study. 
Appl. Math. Comput. [submitted] (2019)
\bibitem{MarkovEvolution}
Finkelshtein, D., Kondratiev, Y., Oliveira, M.J.: Markov evolutions and hierarchical equations 
in the continuum. I. One-component systems. J. Evol. Equ. 9, 197--233 (2009)
\bibitem{VlasovScaling}
Finkelshtein, D., Kondratiev, Y., Kutoviy, O.: Vlasov scaling for stochastic dynamics of continuous 
systems. J. Stat. Phys. 141, 158--178 (2010)
\bibitem{BanasiakLachowicz}
Banasiak, J., Lachowicz, M.: Methods of Small Parameter in Mathematical Biology. 
Birh\"auser/Springer, Cham (2014)
\bibitem{SLM}
Finkelshtein, D., Kondratiev, Yu., Kozitsky, Yu., Kutoviy, O.: The statistical dynamics of a spatial
logistic model and the related kinetic equation. Math. Models Methods Appl. Sci. 25, 343--370 (2015)
\bibitem{PerturbativeApproaches}
Finkelshtein, D., Kondratiev, Y., Kutoviy, O.: Statistical dynamics of continuous systems: perturbative
and approximative approaches. Arabian Journal of Mathematics 4, 255--300 (2015)
\bibitem{SPM}
Kondratiev, Yu., Kozitsky, Yu.: The evolution of states in a spatial pupulation model. J. Dynam. Differ.
Equ. 30, 135--173 (2018)
\bibitem{NumericalMethods}
Chapra, S.C., Canale, R.P.: Numerical Methods for Engineers, 7th edn. McGraw-Hill Education, Penn
Plaza, New York (2015)
\bibitem{BeyondMeanfield}
Omelyan, I., Kozitsky, Yu.: Spatially inhomogeneous population dynamics: beyond the mean field approximation. J. Phys. A.: Math. Theor. 52, 305601 (2019)
\bibitem{BeyondKirkwood}
Omelyan, I.: Spatial population dynamics: beyond the Kirkwood superposition approximation by advancing to the Fisher-Kopeliovich ansatz. ArXiv:1907.00223, Physica A [submitted] (2019)
\bibitem{NumericalRecipes}
Press, W.H., Teukolsky S.A., Vetterling W.T., Flannery B.P.: Numerical Recipes, The Art of Scientific
Computing, 3rd edn., Cambridge University Press, Cambridge (2007)
\bibitem{Convolution2D}
Baddour, N.: Operational and convolution properties of two-dimensional Fourier transforms in polar
coordinates. J. Opt. Soc. Am. A 26, 1767--1777 (2009)
\bibitem{Convolution3D}
Baddour, N.: Operational and convolution properties of three-dimensional Fourier transforms in spherical polar coordinates. J. Opt. Soc. Am. A 27, 2144--2155 (2010)
\bibitem{Nikitin}
Nikitin, A.A., Nikolaev M.V.: Equilibrium integral equations with Kurtosian kernels in spaces of various
dimensions. Mosc. Univ. Comput. Math. Cybern. 42, 105--113 (2018)

\end{thebibliography}


\end{document}